\documentclass{cmslatex}
\usepackage{latexsym, amssymb, enumerate, amsmath}

\usepackage{graphicx}
\usepackage{cite}

\renewcommand{\theequation}{\arabic{section}.\arabic{equation}}

\newcommand{\RR}{\mathbb{R}}
\newcommand{\bs}{\boldsymbol}
\newcommand{\x}{\boldsymbol{x}}
\newcommand{\y}{\boldsymbol{y}}
\newcommand{\z}{\boldsymbol{z}}
\newcommand{\f}{\boldsymbol{f}}
\newcommand{\g}{\boldsymbol{g}}
\newcommand{\bx}{\boldsymbol{\bar x}}
\newcommand{\bz}{\boldsymbol{\bar z}}
\newcommand{\LL}{\boldsymbol{L}}
\newcommand{\Lx}{\boldsymbol{L}_x}
\newcommand{\Ly}{\boldsymbol{L}_y}
\newcommand{\C}{\boldsymbol{C}}
\newcommand{\CC}{\boldsymbol{\mathcal C}}

\newcommand{\F}{\boldsymbol{F}}
\newcommand{\bF}{\boldsymbol{\bar F}}
\newcommand{\G}{\boldsymbol{G}}
\newcommand{\R}{\boldsymbol{R}}
\newcommand{\dt}{\mathrm{d}t}
\newcommand{\dX}{\mathrm{d}\x}
\newcommand{\dY}{\mathrm{d}\y}
\newcommand{\dZ}{\mathrm{d}\z}
\newcommand{\dtau}{\mathrm{d}\tau}

\newcommand{\ds}{\mathrm{d}s}
\newcommand{\mux}{\mu_{\x}}
\newcommand{\dmux}{\mathrm{d}\mux}

\newcommand{\dif}{\mathrm{d}}
\newcommand{\parderiv}[2]{\frac{\partial #1}{\partial #2}}

\begin{document}

\title{The response of reduced models of multiscale dynamics to
small external perturbations
\thanks{
}}
\author{Rafail V. Abramov \thanks {Department of Mathematics,
    Statistics and Computer Science, University of Illinois at
    Chicago, 851 S. Morgan st., Chicago, IL 60607
    (abramov@math.uic.edu).  http://www.math.uic.edu/\~{}abramov} \and
  Marc Kjerland \thanks {Department of Mathematics, Statistics and
    Computer Science, University of Illinois at Chicago, 851 S. Morgan
    st., Chicago, IL 60607 (mkjerl2@uic.edu).
    http://www.math.uic.edu/\~{}kjerland}}



\pagestyle{myheadings} \markboth{The response of reduced models to
small external perturbations}{Rafail V. Abramov and Marc Kjerland}
\maketitle

\begin{abstract}
In real-world geophysical applications (such as predicting the climate
change), the reduced models of real-world complex multiscale dynamics
are used to predict the response of the actual multiscale climate to
changes in various global atmospheric and oceanic parameters. However,
while a reduced model may be adjusted to match a particular dynamical
regime of a multiscale process, it is unclear why it should respond to
external perturbations in the same way as the underlying multiscale
process itself. In the current work, the authors study the statistical
behavior of a reduced model of the linearly coupled multiscale Lorenz
96 system in the vicinity of a chosen dynamical regime by perturbing
the reduced model via a set of forcing parameters and observing the
response of the reduced model to these external perturbations.
Comparisons are made to the response of the underlying multiscale
dynamics to the same set of perturbations. Additionally, practical
viability of linear response approximation via the
Fluctuation-Dissipation theorem is studied for the reduced model.
\end{abstract}


\section{Introduction}\label{intro}

A reduced model for slow variables of multiscale dynamics is a
lower-dimensional dynamical system, which ``resolves'' (that is,
qualitatively approximates in some appropriate sense) major large
scale slow variables of the underlying higher-dimensional multiscale
dynamics while at the same time being relatively simple and
computationally inexpensive to work with. This is important in
real-world applications of contemporary science, such as geophysical
science and climate change prediction
\cite{FraMajVan,Has,BuiMilPal,Pal3,BraBer,MajFraCro,KraKonGhi}, where
the actual underlying physical process is impossible to model
directly, and its reduced approximation has to be designed for such a
purpose. Reduced dynamics were used to model global circulation
patterns \cite{FraMaj,Bra,NewSarPen,WhiSar,ZhaHel}, and large-scale
features of tropical convection \cite{MajKho,KhoMajKat}. Typically,
reduced models of multiscale dynamics consist of simplified
lower-dimensional dynamics of the original multiscale dynamics for the
resolved variables, with additional terms and parameters which serve
as replacements to the missing coupling terms with the unresolved
variables of the underlying physical process. These extra parameters
in the reduced model are usually computed to match a particular
dynamical regime of the underlying multiscale dynamics
\cite{Abr9,Abr10,Abr11,CroVan,FatVan,MajTimVan,MajTimVan2,MajTimVan3,MajTimVan4,Wilks,KatVla,AzeBerTim}.
In particular, if the underlying multiscale process changes its
dynamical regime (for example, in response to changes in its own
forcing parameters), then the parameters of the corresponding reduced
model have to be appropriately readjusted to match its dynamical
regime to the new regime of the multiscale dynamics.

In some real-world applications, such as the climate change
prediction, the reduced models of complex multiscale climate dynamics
are used to predict the response of the actual multiscale climate to
changes in various global atmospheric and oceanic parameters. However,
while a reduced model may be manually adjusted to match a particular
dynamical regime of a multiscale process, it is unclear whether it
should respond to identical external perturbations {\em a priori} in
the same way as the multiscale process, without any extra
readjustments. How do reduced models of multiscale dynamics, adjusted
to a particular dynamical regime, respond to external perturbations
which force them out of this regime?  Is their response similar to the
response of the underlying multiscale dynamics to the same external
perturbation?  It is quite clear that the reduced dynamics evolve on a
set with lower dimension than that of the full multiscale
dynamics. How do the properties of this limiting set respond to
changing external forcing parameter, in comparison to the full
multiscale attractor?

Here we develop a set of criteria for similarity of the response to
small external perturbations between slow variables of multiscale
dynamics and those of a reduced model for slow variables only,
determined through statistical properties of the unperturbed
dynamics. We also carry out a computational study of the difference in
responses of the full multiscale and deterministic reduced dynamics of
the linearly coupled rescaled Lorenz '96 model from \cite{Abr9,Abr11}
to identical external perturbations.  We compare and contrast both the
actual (``ideal'') responses of the multiscale and reduced models
directly to finitely small perturbations of external forcing, and the
linear response predictions of the reduced models via the
Fluctuation-Dissipation theorem
\cite{Abr5,Abr6,Abr7,AbrMaj4,AbrMaj5,AbrMaj6,AbrMaj7,MajAbrGro,Ris}.
Two different types of forcing perturbations are used: the
time-independent Heaviside forcing, and the simple time-dependent ramp
forcing.  The manuscript is structured as follows. In Section
\ref{sec:averaged} we formulate the standard averaging formalism to
obtain the averaged slow dynamics from a general two-scale dynamical
system. Section \ref{sec:response} describes statistically tractable
criteria to ensure similarity of responses between a two-scale system
and its averaged slow dynamics.  In Section \ref{sec:implement} we
describe the first-order reduced model approximation to a two-scale
dynamics with linear coupling between the slow and fast variables,
previously developed in \cite{Abr9,Abr10,Abr11}. In Section
\ref{sec:lorenz96} we introduce the two-scale Lorenz '96 toy model
which will be our testbed for this method.  In Section
\ref{sec:results} we present comparisons of the large features of the
multiscale and reduced systems, including statistical comparisons as
well as the ability of the reduced model to capture perturbation
response of the multiscale system.  Section \ref{sec:summary}
summarizes the results and suggests future work.

\section{Averaged slow dynamics for a general two-scale system}
\label{sec:averaged}

A general two-scale dynamical system with slow variables $\x$ and fast
variables $\y$ is usually represented as
\begin{equation}
    \left\{\begin{aligned}
    &\frac{\dX}{\dt} = \ \F(\x,\y),\\
    &\frac{\dY}{\dt} = \ \G(\x,\y),
    \end{aligned}\right.
    \label{eq:dyn_sys}
\end{equation}
where, $\x(t)\in\RR^{N_x}$ are the slow variables of the system,
$\y(t)\in\RR^{N_y}$ are the fast variables, and $\F$ and $\G$ are
nonlinear differentiable functions. The integer parameters $N_x\ll
N_y$ are the dimensions of the slow and fast variable subspaces,
respectively. Usually, a time-scale separation parameter is used to
denote the difference in time scales between the slow and fast
variables, however, here we omit it, as the framework for reduced
models from \cite{Abr9,Abr10,Abr11}, which we use here, does not
require such a parameter to be explicitly present.

Under the assumption of ``infinitely fast'' $\y$-variables, one can
use the averaging formalism \cite{Pap,Van,Vol,PavStu} to write the
averaged system for slow variables alone:
\begin{equation}
    \frac{\dX}{\dt} = \bF(\x),\qquad\bF(\x)=\int\F(\x,\y)\dmux(\y),
    \label{eq:dyn_sys_slow_averaged}
\end{equation}
where $\mux$ is the invariant distribution measure of the fast
limiting system
\begin{equation}
    \frac{\dZ}{\dtau} = \G(\x,\z),
    \label{eq:dyn_sys_fast_limiting}
\end{equation}
with $\x$ above in \eqref{eq:dyn_sys_fast_limiting} being a constant
parameter. We express the slow solutions of the two-scale system in
\eqref{eq:dyn_sys} and the averaged system in
\eqref{eq:dyn_sys_slow_averaged} in terms of differentiable flows:
\begin{subequations}
\begin{equation}
\x(t)=\phi^t(\x_0,\y_0)\quad\mbox{for the two-scale system},
\end{equation}
\begin{equation}
\x_A(t)=\phi_A^t(\x_0)\quad\mbox{for the averaged system}.
\end{equation}
\end{subequations}
It can be shown (see \cite{Pap,Van,Vol,PavStu} and references therein)
that if the time scale separation between $\x$ and $\y$ is large
enough, then, for the identical initial conditions $\x_0$ and generic
choice of $\y_0$, the solution $\x_A(t)$ of the averaged system in
\eqref{eq:dyn_sys_slow_averaged} remains near the solution $\x(t)$ of
the original two-scale system in \eqref{eq:dyn_sys} for finitely long
time.

\section{Criteria of similarity of responses to small external
perturbations for general two-scale system and its averaged slow
dynamics}
\label{sec:response}

Let $\mu$ and $\mu_A$ denote the invariant distribution measures for
the two-scale system in \eqref{eq:dyn_sys} and the averaged system in
\eqref{eq:dyn_sys_slow_averaged}, respectively. Also, let $h(\x)$ be a
differentiable test function. Then, the statistically average values
of $h$ for both two-scale and averaged systems are given via
\begin{subequations}
\label{eq:h}
\begin{equation}
\langle h\rangle=\int h(\x)\dif\mu(\x,\y),
\end{equation}
\begin{equation}
\langle h\rangle_A=\int h(\x)\dif\mu_A(\x).
\end{equation}
\end{subequations}
Now, consider the two-scale system in \eqref{eq:dyn_sys}, and the
averaged system in \eqref{eq:dyn_sys_slow_averaged}, perturbed at the slow
variables by a small time-dependent forcing $\delta\f(t)$:
\begin{subequations}
\label{perturbed-eqn}
\begin{equation}
    \left\{\begin{aligned}
    &\frac{\dX}{\dt} = \ \F(\x,\y)+\delta\f(t),\\
    &\frac{\dY}{\dt} = \ \G(\x,\y),
    \end{aligned}\right.
\end{equation}
\begin{equation}
    \frac{\dX}{\dt} = \bF(\x)+\delta\f(t).
\end{equation}
\end{subequations}
Then, for small enough $\delta\f(t)$, the average responses
$\delta\langle h\rangle(t)$ and $\delta\langle h\rangle_A(t)$ for the
two-scale system in \eqref{eq:dyn_sys} and the averaged system in
\eqref{eq:dyn_sys_slow_averaged}, respectively, can be approximated by
the following linear response relations:
\begin{subequations}
\begin{equation}
\delta\langle h\rangle(t)=\int_0^t\R(t-s)\delta\f(s)\ds,\qquad
\R(t)=\int\parderiv{h(\phi^t(\x,\y))}{\x}\dif\mu(\x,\y),
\end{equation}
\begin{equation}
\delta\langle h\rangle_A(t)=\int_0^t\R_A(t-s)\delta\f(s)\ds,\qquad
\R_A(t)=\int\parderiv{h(\phi_A^t(\x))}{\x}\dif\mu_A(\x).
\end{equation}
\end{subequations}
For details, see
\cite{Abr5,Abr6,Abr7,AbrMaj4,AbrMaj5,AbrMaj6,AbrMaj7,Rue2}.  Above, it
is clear that any differences between $\delta\langle h\rangle(t)$ and
$\delta\langle h\rangle_R(t)$ are due to differences between $\R(t)$
and $\R_A(t)$, since $\delta\f$ is identical in both cases. The
differences between $\R(t)$ and $\R_A(t)$ are, in turn, caused by the
differences between the flows $\phi^t$ and $\phi_A^t$, and the
differences between the invariant distribution measures $\mu$ and
$\mu_A$, which are difficult to quantify in practice. In what follows
we express the differences between $\R(t)$ and $\R_A(t)$ via
statistically tractable quantities. First, we assume that the
invariant measures $\mu$ and $\mu_A$ are absolutely continuous with
respect to the Lebesgue measure, with distribution densities
$\rho(\x,\y)$ and $\rho_A(\x)$, respectively:
\begin{subequations}
\begin{equation}
\R(t)=\int\parderiv{h(\phi^t(\x,\y))}{\x}\rho(\x,\y)\dif\x\dif\y,
\end{equation}
\begin{equation}
\R_A(t)=\int\parderiv{h(\phi_A^t(\x))}{\x}\rho_A(\x)\dif\x.
\end{equation}
\end{subequations}
While it is known that purely deterministic processes may not have
Lebesgue-continuity of their invariant measures \cite{Rue,Rue1,You},
however, even small amounts of random noise, which is always present
in real-world complex geophysical dynamics, usually ensure the
existence of the distribution density. Integration by parts yields
\begin{subequations}
\label{eq:R}
\begin{equation}
\R(t)=-\int h(\phi^t(\x,\y))\parderiv{\rho(\x,\y)}{\x}\dif\x\dif\y,
\end{equation}
\begin{equation}
\R_A(t)=-\int h(\phi_A^t(\x))\parderiv{\rho_A(\x)}{\x}\dif\x.
\end{equation}
\end{subequations}
At this point, let us express $\rho(\x,\y)$ as the product of its
marginal distribution $\bar\rho(\x)$, defined as
\begin{equation}
\bar\rho(\x)=\int\rho(\x,\y)\dif\y,
\end{equation}
and conditional distribution $\rho(\y|\x)$, given by
\begin{equation}
\rho(\y|\x)=\frac{\rho(\x,\y)}{\bar\rho(\x)}.
\end{equation}
It is easy to check that the conditional distribution $\rho(\y|\x)$ satisfies the identity
\begin{equation}
\label{eq:cond_dist}
\int\rho(\y|\x)\dif\y=1\quad\mbox{for all }\x.
\end{equation}
Now, the formula for the linear response operator $\R(t)$ above can be
written as
\begin{equation}
\R(t)=-\int h(\phi^t(\x,\y))\rho(\y|\x)\parderiv{\bar\rho(\x)}{\x}\dif\y\dif\x-
\int h(\phi^t(\x,\y))\parderiv{\rho(\y|\x)}{\x}\bar\rho(\x)\dif\y\dif\x.
\label{eq:R_temp}
\end{equation}
We now denote
\begin{equation}
\label{eq:epsilon}
\varepsilon^t(\x,\y)=\phi^t(\x,\y)-\phi_A^t(\x),
\end{equation}
where $\varepsilon^t(\x,\y)$ is small compared to either
$\phi^t(\x,\y)$ or $\phi_A^t(\x)$ for relevant values of $t$, $\x$ and
$\y$. Then, for the second integral in the right-hand side of
\eqref{eq:R_temp} we write
\begin{equation}
\begin{split}
-\int h(\phi^t(\x,\y))\parderiv{\rho(\y|\x)}{\x}\bar\rho(\x)\dif\y\dif\x=
-\int \left(\int\parderiv{\rho(\y|\x)}{\x}\dif\y\right)h(\phi_A^t(\x))
\bar\rho(\x)\dif\x-\\-\int \nabla h(\phi_A^t(\x))\varepsilon^t(\x,\y)
\parderiv{\rho(\y|\x)}{\x}\bar\rho(\x)\dif\y\dif\x=O(\|\varepsilon\|),
\end{split}
\end{equation}
where the first integral in the right-hand side is zero due to the
condition in \eqref{eq:cond_dist}. Neglecting the $O(\|\varepsilon\|)$
terms in \eqref{eq:R}, we write
\begin{subequations}
\begin{equation}
\R(t)=-\int h(\phi^t(\x,\y))\rho(\y|\x)\parderiv{\bar\rho(\x)}{\x}\dif\y\dif\x,
\end{equation}
\begin{equation}
\R_A(t)=-\int h(\phi_A^t(\x))\parderiv{\rho_A(\x)}{\x}\dif\x.
\end{equation}
\end{subequations}
At this point, we express $\bar\rho(\x)$ and $\rho_A(\x)$ as exponentials
\begin{equation}
\bar\rho(\x)=e^{-\bar b(\x)},\qquad\rho_A(\x)=e^{-b_A(\x)},
\end{equation}
where $\bar b(\x)$ and $b_A(\x)$ are smooth functions, growing to
infinity as $\x$ becomes infinite. The latter yields
\begin{subequations}
\begin{equation}
\R(t)=\int h(\phi^t(\x,\y))\parderiv{\bar b(\x)}{\x}\rho(\x,\y)\dif\y\dif\x,
\end{equation}
\begin{equation}
\R_A(t)=\int h(\phi_A^t(\x))\parderiv{b_A(\x)}{\x}\rho_A(\x)\dif\x.
\end{equation}
\end{subequations}
Replacing invariant measure averages with long-term time averages
yields the following time correlation functions:
\begin{subequations}
\label{eq:R2}
\begin{equation}
\R(t)=\lim_{r\to\infty}\frac 1r\int_0^r h(\x(s+t))\parderiv{\bar b}{\x}(\x(s))\ds,
\end{equation}
\begin{equation}
\R_A(t)=\lim_{r\to\infty}\frac 1r\int_0^r h(\x_A(s+t))\parderiv{b_A}{\x}(\x_A(s))\ds.
\end{equation}
\end{subequations}
Taking into account the arbitrariness of $h$, we conclude that, in
order for $\R_A(t)$ to approximate $\R(t)$ despite the fact that, for
long times $s$, $\x_A(s)$ diverges from $\x(s)$ even for identical
initial conditions, we generally need three conditions to be
approximately satisfied:
\begin{enumerate}
\item For identical initial conditions, $\x_A(t)$ should approximate
  $\x(t)$ (that is, $\varepsilon^t(\x,\y)$ in \eqref{eq:epsilon}
  should indeed be small) on the finite time scale of decay of the
  correlation functions in \eqref{eq:R2};
\item $b_A(\x)$ should approximate $\bar b(\x)$, which means that the
  invariant distribution $\rho_A(\x)$ of the averaged system in
  \eqref{eq:dyn_sys_slow_averaged} should be similar to the
  $\x$-marginal $\bar\rho(\x)$ of the invariant distribution of the
  two-scale dynamical system in \eqref{eq:dyn_sys};
\item The time autocorrelation functions of the averaged system in
  \eqref{eq:dyn_sys_slow_averaged} should be similar to the time
  autocorrelation functions of the slow variables of the two-scale
  system in \eqref{eq:dyn_sys}.
\end{enumerate}
As a side note, observe that the nature of dependence of the
conditional distribution $\rho(\y|\x)$ on $\x$ does not play any role
in the criteria for the similarity of responses. In particular, the
exact factorization of $\rho(\x,\y)$ into its $\x$- and $\y$-marginals
(which means that $\rho(\y|\x)$ is independent of $\x$) is not
required, unlike what was suggested in \cite{MajGerYua} for the
Gaussian invariant states.

\section{Practical implementation of the reduced model
for a two-scale process with linear coupling}
\label{sec:implement}

As formulated above in Sections \ref{sec:averaged} and
\ref{sec:response}, the criteria of the response similarity are
applicable for a broad range of dynamical systems with general forms
of coupling and their averaged slow dynamics. However, the practical
computation of the reduced model approximation to averaged slow
dynamics depends on the form of coupling in the two-scale system
\cite{Abr9,Abr10,Abr11}. In this work, we consider the linear coupling
between the slow and fast variables in the two-scale system
\eqref{eq:dyn_sys}. The linear coupling is the most basic form of
coupling in physical processes, however, because of that it is also
probably the most common form of coupling. For the linear coupling,
the reduced model is constructed according to the method developed
previously in \cite{Abr9}, which we briefly sketch below.

We consider the special setting of \eqref{eq:dyn_sys} with linear
coupling between $\x$ and $\y$:
\begin{equation}
    \left\{\begin{aligned}
    &\frac{\dX}{\dt} = \ \f(\x)+\Ly\y,\\
    &\frac{\dY}{\dt} = \ \g(\y)+\Lx\x,
    \end{aligned}\right.
    \label{multiscaleODE}
\end{equation}
where $\f$ and $\g$ are nonlinear differentiable functions, and $\Lx$
and $\Ly$ are constant matrices of appropriate sizes. The
corresponding averaged dynamics for slow variables from
\eqref{eq:dyn_sys_slow_averaged} simplifies to
\begin{equation}
    \frac{\dX}{\dt} = \f(\x)+\Ly\bz(\x),
    \label{averagedODE}
\end{equation}
where $\bz(\x)$ is the statistical mean state of the fast limiting
system
\begin{equation}
    \frac{\dZ}{\dtau} = \g(\z)+\Lx\x,
    \label{fastlimitingsystem}
\end{equation}
with $\x$ treated as constant parameter. In general, the exact
dependence of $\bz(\x)$ on $\x$ is unknown, except for a few special
cases like the Ornstein-Uhlenbeck process \cite{OrnUhl}. Here, like in
\cite{Abr9,Abr11}, we approximate $\bz(\x)$ via the linear expansion
\begin{equation}
\bz(\x)\approx\bz^*+\CC\Lx(\x-\x^*),
\end{equation}
where $\x^*$ is the statistical average state of the full multiscale
system in \eqref{multiscaleODE}, and $\bz^*=\bz(\x^*)$. The constant
matrix $\CC$ is computed as the time integral of
the correlation function
\begin{equation}
\CC=\left(\int_0^t\C(s)\ds\right)\C^{-1}(0),
\qquad\C(s)=\lim_{r\to\infty}\frac 1r\int_0^r\z(t+s)\z^T(t)\dt,
\label{corrfunc}
\end{equation}
where $\z(t)$ is the solution of \eqref{fastlimitingsystem} for
$\x=\x^*$. The above formula constitutes the quasi-Gaussian
approximation to the linear response of $\bz$ to small constant
forcing perturbations in \eqref{fastlimitingsystem}, and is a good
approximation when the dynamics in \eqref{fastlimitingsystem} are
strongly chaotic and rapidly mixing
\cite{Abr5,Abr6,Abr7,AbrMaj4,AbrMaj5,AbrMaj6,AbrMaj7,MajAbrGro,Ris}.
With \eqref{corrfunc}, the reduced system in \eqref{averagedODE}
becomes the explicitly defined deterministic reduced model for slow
variables alone:
\begin{equation}
    \frac{\dX}{\dt} = \f(\x)+\Ly\bz^*+\LL(\x-\x^*),
    \label{reducedODE}
\end{equation}
where $\LL=\Ly\CC\Lx$. In what follows, the ``zero-order'' model
refers to \eqref{reducedODE} with the last term set to zero (such that
the coupling is parameterized only by the constant term
$\Ly\bz^*$). For details, see \cite{Abr9,Abr11} and references
therein.

\section{Testbed -- the rescaled Lorenz '96 system}
\label{sec:lorenz96}

In the current work, we test the response of the reduced model for
slow variables on the rescaled Lorenz '96 system with linear coupling
\cite{Abr9}, which is obtained from the original two-scale Lorenz '96
system \cite{Lor} by appropriately rescaling dynamical variables to
approximately set their mean states and variances to zero and one,
respectively. Below we present a brief exposition of how the rescaled
Lorenz '96 model is derived.

\subsection{The original two-scale Lorenz '96 system}

The original two-scale Lorenz '96 system \cite{Lor} is given by
\begin{equation}
\left\{\begin{aligned} &\dot x_i = x_{i-1}(x_{i+1}-x_{i-2}) - x_i +
F_x - \frac{\lambda_y}{J}\sum_{j=1}^J y_{i,j},\\ &\dot y_{i,j} =
\frac{1}{\varepsilon}\left[y_{i,j+1}(y_{i,j-1}-y_{i,j+2}) - y_{i,j} +
  F_y + \lambda_x x_i\right],
\end{aligned}\right.
\end{equation}
where $1\leq i\leq N_x, 1\leq j\leq J,$ and periodic boundary
conditions $x_{i+N_x}=x_i$, $y_{i+N_x,j}=y_{i,j}$ and
$y_{i,j+J}=y_{i+1,j}.$ Here $F_x$ and $F_y$ are constant forcing
terms, $\lambda_x$ and $\lambda_y$ constant coupling parameters, and
$\varepsilon$ is the time scale separation parameter. Throughout this
paper we will consider systems with twenty slow variables $(N_x=20)$
and eighty fast variables $(N_y=80,\ J=4)$.

In Lorenz's original formulation \cite{Lor} studying predictability in
atmospheric-type systems, he begins with the uncoupled system
\begin{equation}
    \dot x_i = x_{i-1}(x_{i+1}-x_{i-2}) - x_i + F, \qquad
    i=1,\ldots,N,
\end{equation}
with periodic boundary conditions. This system has generic features of
geophysical flows, namely a nonlinear advection-like term, linearly
unstable waves, damping, forcing, mixing, and chaos \cite{MajWan}. The
simple formulation, with invariance under index translation and a
uniform forcing term $F$, allows for straightforward analysis - in
particular the long-time statistics of each variable should be
identical and will only depend on $F$. Additionally, the chaos and
mixing of the system are simply regulated by the forcing, with
decaying solutions for $F$ near zero, periodic solutions for $F$
slightly larger, weakly chaotic quasi-periodic solutions around $F=6$,
and chaotic and strongly mixing systems around $F=16$ and
higher. Lorenz's two-time coupled system was introduced to study
predictability and Lyapunov exponents of systems with subgrid
phenomena on faster timescales, and one of the authors of the current
work has recent results showing that coupling two chaotic systems can
suppress chaos in the slower system \cite{Abr8}.

\subsection{Rescaled Lorenz '96 model}

To simplify the analysis of coupling trends for the two-time system,
we will scale out the dependence of the mean state and mean energy on
the forcing term $F$. Due to the translational invariance, the
long-term mean $\overline{x}$ and standard deviation $\sigma$ for the
uncoupled system are the same for all $x_i$. So we rescale $\x$ and
$t$ as
\[
x_i= \overline{x} + \sigma \hat{x}_i, \quad t=\frac{\tau}{\sigma},
\]
where the new variables $\hat{x}$ have zero mean and unit standard
deviation, while their time autocorrelation functions have normalized
scaling across different dynamical regimes (that is, different
forcings $F$) for short correlation times. This rescaling was
previously used in \cite{MajAbrGro}. In the rescaled variables, the
uncoupled Lorenz model becomes
\begin{equation}
\dot{\hat{x}}_i =
\left(\hat{x}_{i-1}+\frac{\overline{x}}{\sigma}\right)(\hat{x}_{i+1}-\hat{x}_{i-2})
- \frac{\hat{x}_i}{\sigma} + \frac{F-\overline{x}}{\sigma^2},
\end{equation}
where $\overline{x}$ and $\sigma$ are functions of $F$.

We similarly rescale the coupled two-scale Lorenz '96 model:
\begin{equation}
\left\{\begin{aligned} \frac{\mathrm{d}x_i}{\dt} &=
\left(x_{i-1}+\frac{\overline{x}}{\sigma_x}\right)(x_{i+1}-x_{i-2}) -
\frac{x_i}{\sigma_x} + \frac{F_x-\overline{x}}{{\sigma_x}^2} -
\frac{\lambda_y}{J}\sum_{j=1}^J
y_{i,j},\\ \varepsilon\frac{\mathrm{d}y_{i,j}}{\dt} &=
\left(y_{i,j+1}+\frac{\overline{y}}{\sigma_y}\right)(y_{i,j-1}-y_{i,j+2})
- \frac{y_{i,j}}{\sigma_y} + \frac{F_y-\overline{y}}{{\sigma_y}^2} +
\lambda_x x_i,
\end{aligned}\right.
\label{rescaled-l96}
\end{equation}
where $\{\overline{x},\ \sigma_x\}$ and $\{\overline{y},\ \sigma_y\}$
are the long term means and standard deviations of the uncoupled
systems with $F_x$ or $F_y$ as constant forcing, respectively. It is
this rescaled coupled Lorenz '96 system that we focus on for the
closure approximation.

Before any numerical tests, one can already anticipate that the
zero-order reduced system will be inadequate for this model even with
such simple coupling. Once the reference state $\x^*$ is determined and $\overline{\z}^*$
computed, the zero-order reduced system is given, according to
\eqref{reducedODE}, by
\begin{equation}
  \dot{\hat{x}}_i =
  \left(\hat{x}_{i-1}+\frac{\overline{x}}{\sigma}\right)(\hat{x}_{i+1}-\hat{x}_{i-2})
  - \frac{\hat{x}_i}{\sigma} + \frac{F-\overline{x}}{\sigma^2} -
  \lambda_y \overline{z}^*.
\end{equation}
This is equivalent to perturbing $F_x$ by $-\sigma_x^2 \lambda_y
\overline{z}^*$, which we expect to be small since $\hat{x}$ and
$\hat{y}$ have zero mean in the uncoupled setting. In particular, we
expect this perturbation to have only a small effect on the
dynamics. However, in the multiscale dynamics it has been shown that a
chaotic regime in the fast system can suppress chaos when coupled to
the slow system \cite{Abr8}, and this phenomenon is completely lost in
the zero-order model.

\section{Numerical experiments}
\label{sec:results}

Here we compare numerical results of the rescaled two-scale Lorenz '96
system with its corresponding reduced system. In particular we look at
the ability of the reduced system captures some statistical quantities and
how well it captures mean response to perturbations in the slow
variables.

In all parameter regimes considered, we have a slow system consisting
of twenty variables $(N_x=20)$ coupled with a fast system of eighty
variables $(N_y=80)$. We use a fourth order Runge-Kutta method with
timestep $dt=\varepsilon/10$ in the multiscale system and $dt=1/10$ in
the reduced system. To compute the mean response, an ensemble of
$10^4$ points is sampled from a single trajectory which has been
allowed to settle onto the attractor. Using the translational symmetry
of the Lorenz '96 system, we rotate the indices to generate an
ensemble twenty times larger.

On a modern laptop, the initial calculation to generate the reduced
system for the Lorenz '96 system takes only a few minutes; once
computed, numerical simulation of the reduced system is faster than
the multiscale system by a factor on the order of
$\varepsilon^{-1}$. Computing the mean response for a single forcing
for 5 time units with a sufficiently large ensemble size ($10^5$
trajectories) takes over an hour in the multiscale system with
$\varepsilon=10^{-2}$ but less than three minutes for the
corresponding reduced system.

\subsection{Comparison of statistical properties of the two-scale and
reduced systems}
\label{statistics_subsection}

In Section \ref{sec:response} we outlined the main requirements for
correctly capturing the response of the two-scale system by its
reduced model. Those were the approximation of joint distribution
density functions (DDF) for slow variables, and the time
autocorrelation functions of the time series. It is, of course, not
computationally feasible to directly compare the 20-dimensional DDFs
and time autocorrelations for all possible test functions. However, it
is possible to compare the one-dimensional marginal DDFs and simple
time autocorrelations for individual slow variables, to have a rough
estimate on how the statistical properties of the multiscale dynamics
are reproduced by the reduced model.

\begin{figure}[!hbtp]
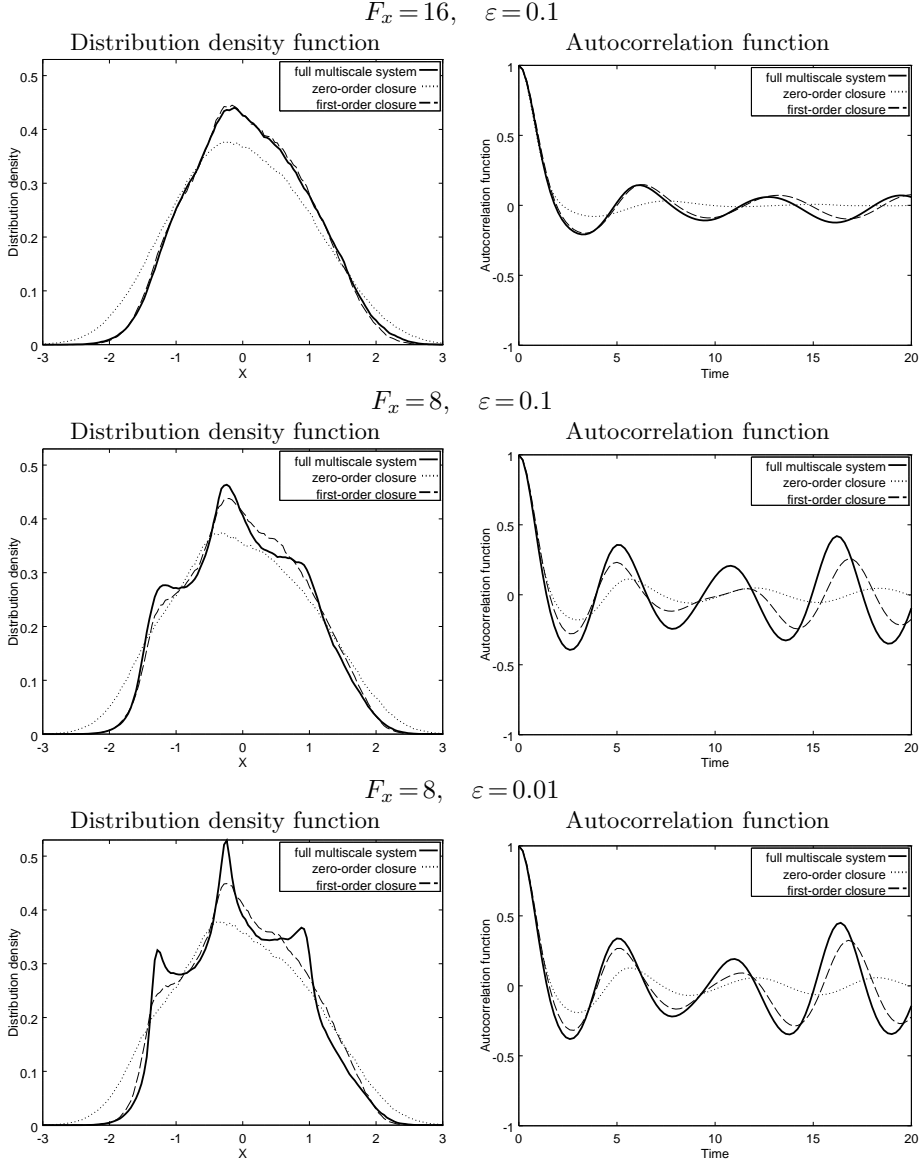

    \center

    $F_x=16,\quad\varepsilon=0.1$\\ 
    \begin{tabular}{cc}
        \small{Distribution density function} & \small{Autocorrelation function}\\
        \includegraphics[width=0.45\textwidth]{{{figures/PDFs_Nx20_Ny80_Fx16_Fy12_lx0.4_ly0.4_eps0.1}}}&
        \includegraphics[width=0.45\textwidth]{{{figures/autocorr_Nx20_Ny80_Fx16_Fy12_lx0.4_ly0.4_eps0.1}}}
    \end{tabular} \\

    $F_x=8,\quad\varepsilon=0.1$\\ 
    \begin{tabular}{cc}
        \small{Distribution density function} & \small{Autocorrelation function}\\
        \includegraphics[width=0.45\textwidth]{{{figures/PDFs_Nx20_Ny80_Fx8_Fy12_lx0.4_ly0.4_eps0.1}}}&
        \includegraphics[width=0.45\textwidth]{{{figures/autocorr_Nx20_Ny80_Fx8_Fy12_lx0.4_ly0.4_eps0.1}}}
    \end{tabular} \\

    $F_x=8,\quad\varepsilon=0.01$\\ 
    \begin{tabular}{cc}
        \small{Distribution density function} & \small{Autocorrelation function}\\
        \includegraphics[width=0.45\textwidth]{{{figures/PDFs_Nx20_Ny80_Fx8_Fy12_lx0.4_ly0.4_eps0.01}}}&
        \includegraphics[width=0.45\textwidth]{{{figures/autocorr_Nx20_Ny80_Fx8_Fy12_lx0.4_ly0.4_eps0.01}}}
    \end{tabular} \\

    \caption{Distribution density and autocorrelation functions of slow variables.}
    \label{DDFs}
\end{figure}

In figure \ref{DDFs} we compare the distribution density functions and
autocorrelation functions of the slow variables. The DDFs are computed
using bin-counting, and the autocorrelation function $\langle
x_i(t)x_i(t+s)\rangle$, averaged over $t$, is normalized by the
variance $\langle x_i^2\rangle$. Results from three parameter regimes
are presented, and in all three regimes the fast system is chaotic and
weakly mixing $(F_y=12)$ and the coupling strength is chosen to be
large enough $(\lambda_x=\lambda_y=0.4)$ so that the multiscale
dynamics are challenging to approximate. Of particular interest are
timescale separations of $\varepsilon=10^{-1}$ and
$\varepsilon=10^{-2}$.

First we consider a chaotic and strongly mixing slow regime
$(F_x=16)$. Figures are presented for the timescale separation
$\varepsilon=10^{-1}$ only, because in this regime the picture is very
similar for $\varepsilon=10^{-2}$. We also consider a weakly chaotic
and quasi-periodic slow regime $(F_x=8)$. In this regime, the coupled
dynamics are more dependent on the timescale separation so we present
results for both $\varepsilon=10^{-1}$ and $\varepsilon=10^{-2}$.
Statistical quantities of other regimes, including regimes with more
periodic behavior, have been presented in \cite{Abr9}.

\begin{figure}
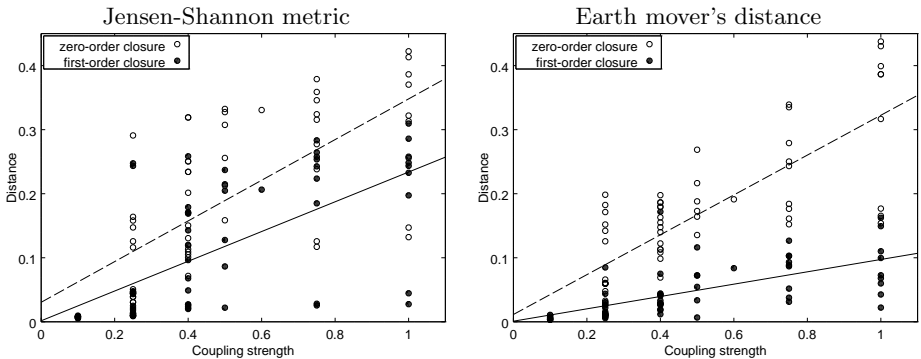

    \begin{tabular}{cc}
        \small{Jensen-Shannon metric} & \small{Earth mover's distance}\\
        \includegraphics[width=0.45\textwidth]{{{figures/JSD}}}&
        \includegraphics[width=0.45\textwidth]{{{figures/EMD}}}
    \end{tabular}
    \caption{Distances between DDFs of reduced and multiscale systems}
    \label{DDF_metrics}
\end{figure}

To more systematically compare DDFs for many parameter regimes, we
introduce two metrics on the space of distributions. First is the
Jensen-Shannon metric which is derived from the information-theoretic
Kullback-Leibler divergence \cite{KulLei} and is given by
\begin{equation}
  m_{\mathrm{JS}}(P,Q) = \frac{1}{\sqrt{2}}\left(\int_{-\infty}^\infty
  \log\left(\frac{2p(x)}{p(x)+q(x)}\right)p(x)dx +
  \int_{-\infty}^\infty
  \log\left(\frac{2q(x)}{p(x)+q(x)}\right)q(x)dx\right)^{1/2},
\end{equation}
where $p$ and $q$ are densities on distributions $P$ and $Q$. The next
metric we consider is the earth mover's distance \cite{RubTomGui},
which measures the minimum energy needed to move one DDF to another as
though they were piles of dirt; the energy cost is the amount of
`dirt' times the distance it moved. One nice property is that the
distance between two delta distributions $\delta_\alpha$ and
$\delta_\beta$ is simply $|\alpha-\beta|$. For distribution functions
of one-dimensional random variables, the earth-mover's distance is the
$L^1$ norm of the difference of the cumulative distributions:
\begin{equation}
    m_{\mathrm{EM}}(P,Q) = \int_{-\infty}^\infty \left|
    \int_{-\infty}^x p(s)-q(s)\mathrm{d}s \right| \mathrm{d}x.
\end{equation}
Figure \ref{DDF_metrics} shows distances between reduced systems DDFs
and the corresponding multiscale slow variable DDFs. A variety of
regimes is considered, with coupling parameters
$\lambda_x,\lambda_y\in[0.1,1]$, forcing parameters
$F_x\in\{6,7,8,10,16\}$ and $F_y\in\{8,12,16\}$, and timescale
separations $\varepsilon\in\{10^{-1},10^{-2}\}$. The data points are
plotted with respect to coupling parameter $\lambda_x$. For each
regime considered, the corresponding distances are shown for both
zero-order and first-order reduced systems.
As the coupling strength between fast and slow systems increases, it
is apparently more difficult for the reduced systems to capture the
correct slow dynamics of the multiscale system, as suggested by the
linear best-fit curves. It should be noted that this correlation is
slightly weaker when plotted against $\lambda_y$, the coupling
parameter for the fast system. However, the distribution densities of
the first-order reduced system are consistently closer in both metrics
than the zero-order system to the multiscale system.

\subsection{Mean state response to small perturbations}
\label{meanresponse_subsection}

\begin{figure}[!hbtp]
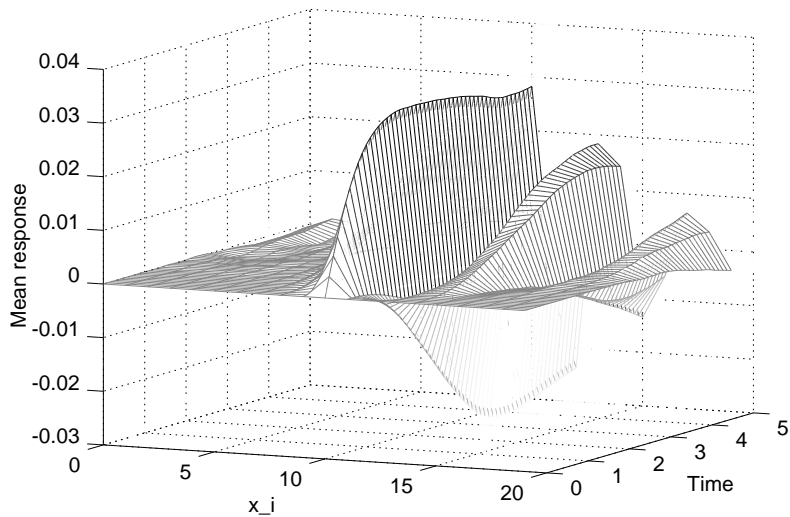

    \center
    \includegraphics[width=0.8\textwidth]{{{figures/meanresponse}}}
    \caption{Mean response of $\x$ to step forcing at node 
      $x_{11}$. $N_x=20,\ N_y=80,\ F_x=16,\ F_y=12,\ \lambda_x=0.4,\ \lambda_y=0.4,\ \varepsilon=0.1$ }
    \label{meanresp_fig}
\end{figure}

In this section we examine the response of the mean state
$\langle\x\rangle$ of the slow variables (that is, $h(\x)=\x$ in
\eqref{eq:h}) in the Lorenz '96 system to two simple types of small
external forcing:
\begin{enumerate}
\item Heaviside step forcing
\begin{equation*}
    \delta\f_{\mathrm{H}}(t) = \begin{cases}
    0& \text{if } t < 0,\\
    \bs{v}_{\mathrm{H}} & \text{if } t > 0,
    \end{cases}
    \label{heaviside}
\end{equation*}
\item Ramp forcing
\begin{equation*}
    \delta\f_{\mathrm{ramp}}(t) = \begin{cases}
    0& \text{if } t < 0,\\
    t \bs{v}_{\mathrm{ramp}} & \text{if } t > 0,
    \end{cases}
    \label{ramp}
\end{equation*}
\end{enumerate}
where $\bs{v}_{\mathrm{H}}$ and $\bs{v}_{\mathrm{ramp}}$ are constant
vectors. To compute the response of the mean state $\langle\x\rangle$,
we generate an initial ensemble sampled from a trajectory that has
been given sufficient time to settle onto the attractor. For each
ensemble member we let evolve a short trajectory under the unperturbed
dynamics as well as the under the perturbed dynamics, then we take the
difference between these two trajectories and average over the entire
ensemble. Here we consider forcing of the form
$\bs{v}=\alpha\hat{e}_j$, where $\alpha$ is constant and $\hat{e}_j$ a
standard basis vector in $\RR^{N_x}$. In the translation-invariant
Lorenz '96 system, without loss of generality we need only consider a
single such vector, say $\hat{e}_0$. Figure \ref{meanresp_fig} shows
the mean response of the slow variables to small Heaviside forcing in
the two-time rescaled Lorenz '96 system, where the magnitude of the
forcing $|\bs{v}_{\mathrm{H}}|$ is $1\%$ of $|\f|$ averaged over the
invariant distribution.

\begin{figure}[!hbtp]
    \center
    \small{Mean relative error of ideal response}
    \begin{tabular}{cc}
    \includegraphics[width=0.45\textwidth]{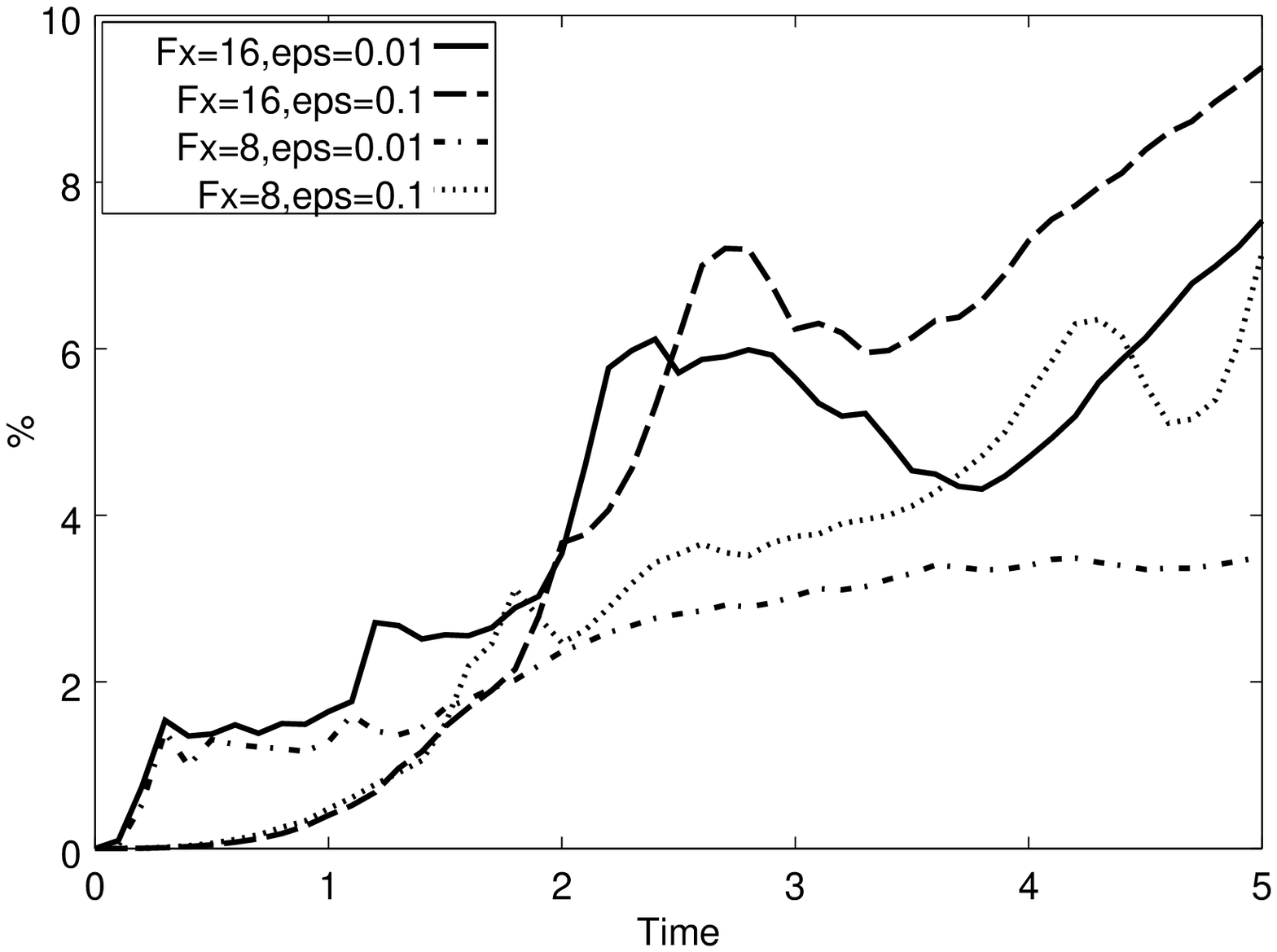} & 
    \includegraphics[width=0.45\textwidth]{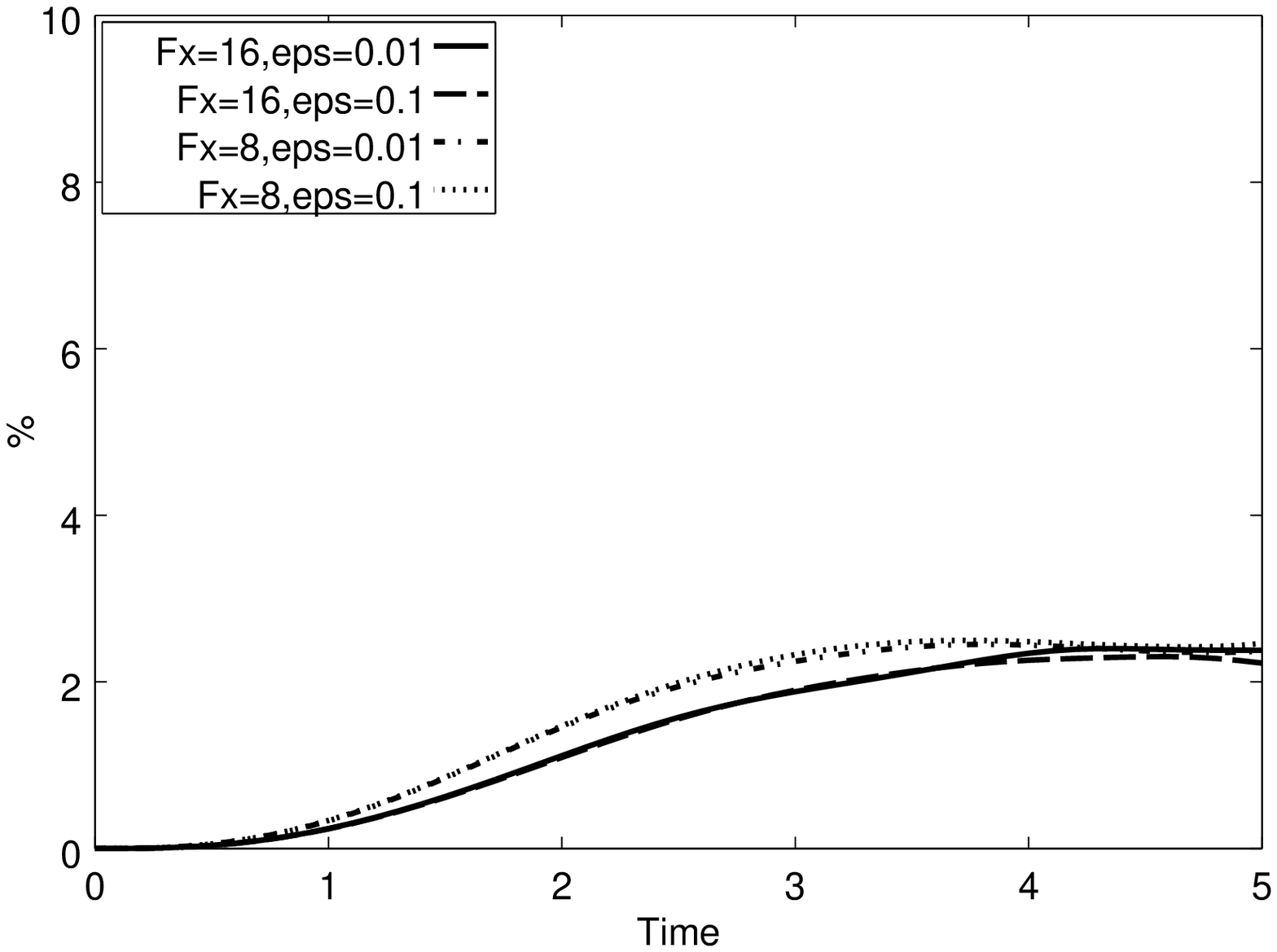} \\
    Multiscale system & Reduced system
    \end{tabular}
    \caption{Nonlinearity of response: relative error ideal response vs actual nonlinear response.}
    \label{nonlinearity}
\end{figure}

Since the external forcing is sufficiently small to consider the
response of the mean state approximately linear, we use the
\emph{ideal} response operator of Gritsun and Dymnikov \cite{GriDym}
(also see
\cite{Abr5,Abr6,Abr7,AbrMaj4,AbrMaj5,AbrMaj6,AbrMaj7,MajAbrGro}) by
generating mean responses for several perturbations and computing the
linear best least-squares fit. This is a time-dependent matrix, and
due to the symmetry of the Lorenz '96 system the dimensionality is
reduced by one so that we have simply a time-dependent vector. With
the ideal response operator, the response to a multitude of forcings
can be readily estimated. We verify the nonlinearity in the actual
response in Figure \ref{nonlinearity}, which shows the growth in time
of the relative error between the ideal response and the actual
response to small Heaviside step forcing, averaged over several
different forcings. We limit the plot to 5 time units because the
large features of the Heaviside response are mostly fully developed by that
time. Note in particular that the response is more linear in the
reduced system, which is probably due to the fact that the Lyapunov
exponents in the reduced system are much smaller than those in the
two-scale system.

\begin{figure}
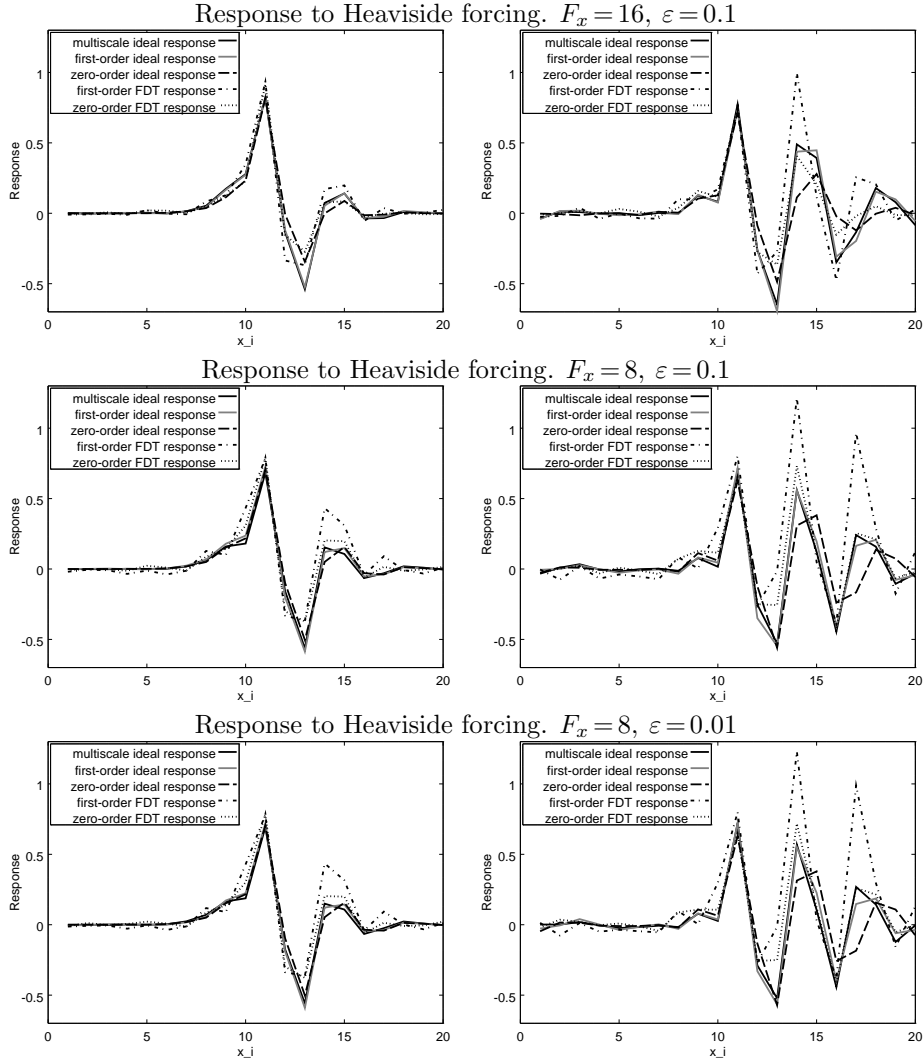

  \center
    Response to Heaviside forcing. $F_x=16,\ \varepsilon=0.1$
    \begin{tabular}{cc}
        \includegraphics[width=0.45\textwidth]{{{figures/response_snapshots_t2_Nx20_Ny80_Fx16_Fy12_lx0.4_ly0.4_eps0.1}}}&
        \includegraphics[width=0.45\textwidth]{{{figures/response_snapshots_t5_Nx20_Ny80_Fx16_Fy12_lx0.4_ly0.4_eps0.1}}}
    \end{tabular}
    Response to Heaviside forcing. $F_x=8,\ \varepsilon=0.1$ \\
    \begin{tabular}{cc}
        \includegraphics[width=0.45\textwidth]{{{figures/response_snapshots_t2_Nx20_Ny80_Fx8_Fy12_lx0.4_ly0.4_eps0.1}}}&
        \includegraphics[width=0.45\textwidth]{{{figures/response_snapshots_t5_Nx20_Ny80_Fx8_Fy12_lx0.4_ly0.4_eps0.1}}}
    \end{tabular}
    Response to Heaviside forcing. $F_x=8,\ \varepsilon=0.01$ \\
    \begin{tabular}{cc}
        \includegraphics[width=0.45\textwidth]{{{figures/response_snapshots_t2_Nx20_Ny80_Fx8_Fy12_lx0.4_ly0.4_eps0.01}}}&
        \includegraphics[width=0.45\textwidth]{{{figures/response_snapshots_t5_Nx20_Ny80_Fx8_Fy12_lx0.4_ly0.4_eps0.01}}}
    \end{tabular}
        \caption{Snapshots of the response operators for the response time
      $T=2$ (left), and $T=5$ (right), Heaviside forcing.}
    \label{fdt_ideal}
\end{figure}

\begin{figure}
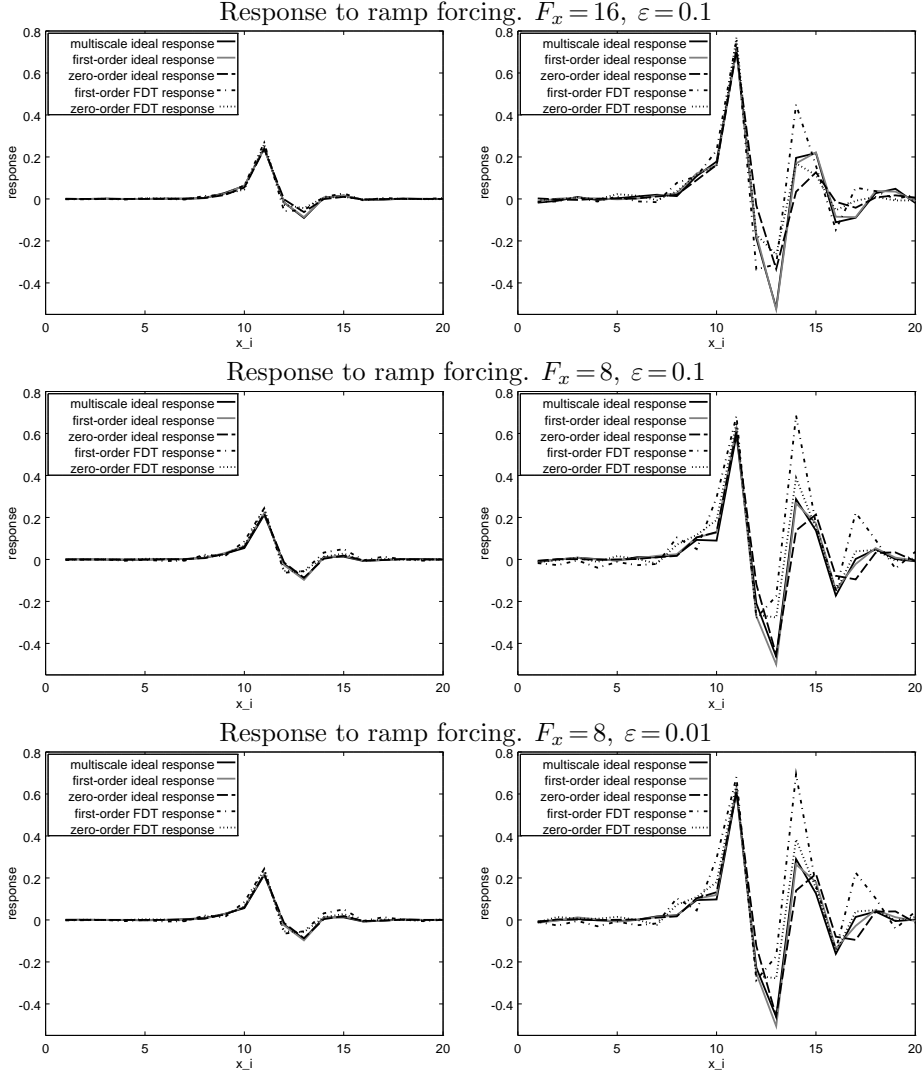

  \center 
    Response to ramp forcing. $F_x=16,\ \varepsilon=0.1$
    \begin{tabular}{cc}
        \includegraphics[width=0.45\textwidth]{{{figures/response_snapshots_ramp_t2_Nx20_Ny80_Fx16_Fy12_lx0.4_ly0.4_eps0.1}}}&
        \includegraphics[width=0.45\textwidth]{{{figures/response_snapshots_ramp_t5_Nx20_Ny80_Fx16_Fy12_lx0.4_ly0.4_eps0.1}}}
    \end{tabular}
    Response to ramp forcing. $F_x=8,\ \varepsilon=0.1$
    \begin{tabular}{cc}
        \includegraphics[width=0.45\textwidth]{{{figures/response_snapshots_ramp_t2_Nx20_Ny80_Fx8_Fy12_lx0.4_ly0.4_eps0.1}}}&
        \includegraphics[width=0.45\textwidth]{{{figures/response_snapshots_ramp_t5_Nx20_Ny80_Fx8_Fy12_lx0.4_ly0.4_eps0.1}}}
    \end{tabular}
    Response to ramp forcing. $F_x=8,\ \varepsilon=0.01$
    \begin{tabular}{cc}
        \includegraphics[width=0.45\textwidth]{{{figures/response_snapshots_ramp_t2_Nx20_Ny80_Fx8_Fy12_lx0.4_ly0.4_eps0.01}}}&
        \includegraphics[width=0.45\textwidth]{{{figures/response_snapshots_ramp_t5_Nx20_Ny80_Fx8_Fy12_lx0.4_ly0.4_eps0.01}}}
    \end{tabular}
    \caption{Snapshots of the response operators for the response time
      $T=2$ (left), and $T=5$ (right), ramp forcing.}
    \label{fdt_ideal_ramp}
\end{figure}

We now compare the ideal responses of the full and reduced
systems. The snapshots of the ideal responses for the two-scale and
reduced models (as well as the linear response approximation,
described in the next section) at times $T=2$ and $T=5$ are shown in
Figures \ref{fdt_ideal} and \ref{fdt_ideal_ramp} for the Heaviside and
ramp forcing, respectively. The response is captured accurately at the
node which is directly forced (here $x_{11}$), but capturing the
off-diagonal response as the perturbation propagates through the
system is more difficult. Indeed, it seems that in the zero-order
reduced systems the propagation speed is slightly faster than in the
two-scale systems, but for the first-order reduced systems the
response is well captured.

\begin{figure}[!hbtp]
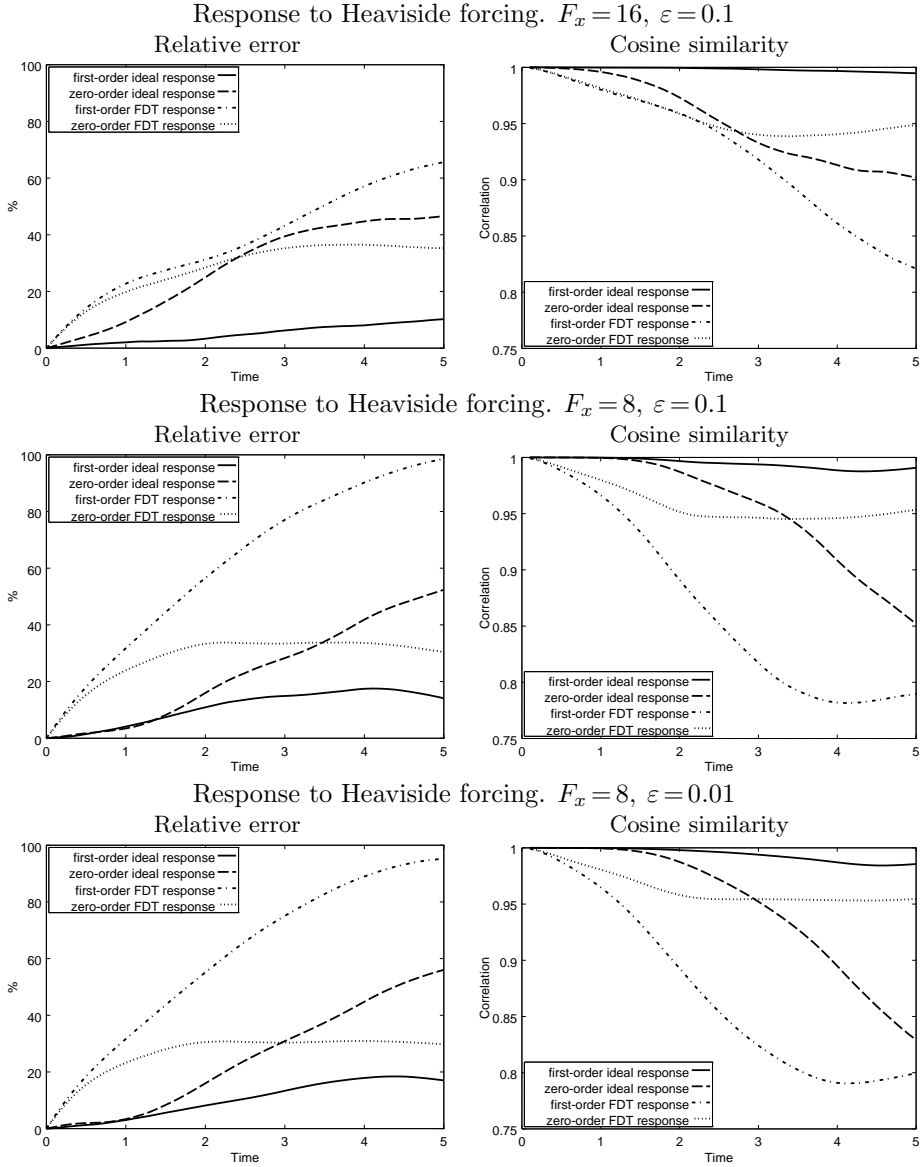

    \center
    Response to Heaviside forcing. $F_x=16,\ \varepsilon=0.1$ 
    \begin{tabular}{cc}
        \small Relative error & \small Cosine similarity \\
        \includegraphics[width=0.45\textwidth]{{{figures/response_error_Nx20_Ny80_Fx16_Fy12_lx0.4_ly0.4_eps0.1}}}&
        \includegraphics[width=0.45\textwidth]{{{figures/response_correlation_Nx20_Ny80_Fx16_Fy12_lx0.4_ly0.4_eps0.1}}}
    \end{tabular}

    Response to Heaviside forcing. $F_x=8,\ \varepsilon=0.1$ 
    \begin{tabular}{cc}
        \small Relative error & \small Cosine similarity \\
        \includegraphics[width=0.45\textwidth]{{{figures/response_error_Nx20_Ny80_Fx8_Fy12_lx0.4_ly0.4_eps0.1}}}&
        \includegraphics[width=0.45\textwidth]{{{figures/response_correlation_Nx20_Ny80_Fx8_Fy12_lx0.4_ly0.4_eps0.1}}}
    \end{tabular}

    Response to Heaviside forcing. $F_x=8,\ \varepsilon=0.01$ 
    \begin{tabular}{cc}
        \small Relative error & \small Cosine similarity\\
        \includegraphics[width=0.45\textwidth]{{{figures/response_error_Nx20_Ny80_Fx8_Fy12_lx0.4_ly0.4_eps0.01}}}&
        \includegraphics[width=0.45\textwidth]{{{figures/response_correlation_Nx20_Ny80_Fx8_Fy12_lx0.4_ly0.4_eps0.01}}}
    \end{tabular}

    \caption{Comparing multiscale ideal response with reduced system ideal \& quasi-Gaussian response operators for Heaviside forcing}
    \label{fdt_ideal_error}
\end{figure}

It is more clear to see the quantitative differences in Figures \ref{fdt_ideal_error}
and \ref{fdt_ideal_ramp_error} which show the
relative distance between the responses as well as their cosine
similarity $\frac{\bs{u}\cdot\bs{v}}{\|\bs{u}\|\|\bs{v}\|}$ versus
time for Heaviside step forcing and ramp forcing, respectively.
Also shown in these figures are
the linear responses computed using the reduced system statistics, as
described in the next section.

\begin{figure}[!hbtp]
    \center
    Response to ramp forcing. $F_x=16,\ \varepsilon=0.1$ 
    \begin{tabular}{cc}
        \small Relative error & \small Cosine similarity\\
        \includegraphics[width=0.45\textwidth]{{{figures/response_error_ramp_Nx20_Ny80_Fx16_Fy12_lx0.4_ly0.4_eps0.1_ramptime5}}}&
        \includegraphics[width=0.45\textwidth]{{{figures/response_correlation_ramp_Nx20_Ny80_Fx16_Fy12_lx0.4_ly0.4_eps0.1_ramptime5}}}
    \end{tabular}

    Ramp forcing. $F_x=8,\ \varepsilon=0.1$ 
    \begin{tabular}{cc}
        \small Relative error & \small Cosine similarity\\
        \includegraphics[width=0.45\textwidth]{{{figures/response_error_ramp_Nx20_Ny80_Fx8_Fy12_lx0.4_ly0.4_eps0.1_ramptime5}}}&
        \includegraphics[width=0.45\textwidth]{{{figures/response_correlation_ramp_Nx20_Ny80_Fx8_Fy12_lx0.4_ly0.4_eps0.1_ramptime5}}}
    \end{tabular}

    Response to ramp forcing. $F_x=8,\ \varepsilon=0.01$ 
    \begin{tabular}{cc}
        \small Relative error & \small Cosine similarity\\
        \includegraphics[width=0.45\textwidth]{{{figures/response_error_ramp_Nx20_Ny80_Fx8_Fy12_lx0.4_ly0.4_eps0.01_ramptime5}}}&
        \includegraphics[width=0.45\textwidth]{{{figures/response_correlation_ramp_Nx20_Ny80_Fx8_Fy12_lx0.4_ly0.4_eps0.01_ramptime5}}}
    \end{tabular}

    \caption{Comparing multiscale ideal response with reduced system ideal \& quasi-Gaussian response operators for ramp forcing}
    \label{fdt_ideal_ramp_error}
\end{figure}
We observe that the first-order reduced system ideal response is a much
closer approximation to the multiscale ideal response than the
corresponding zero-order ideal response. In these regimes the relative
error of the first-order response is limited to about $20\%$ for the
Heaviside forcing and less for the ramp forcing, while in the zero-order
system the error is around $40\%$ for the step forcing and $30\%$ for ramp forcing at time $t=5$.
Remark that in the third plot for ramp forcing response in Figure \ref{fdt_ideal_ramp_error}
there is a small bump in the relative error shortly after the onset of forcing.
This plot corresponds to a weakly chaotic regime ($F_x=8,F_y=12$) with a
large timescale separation ($\varepsilon=10^{-2}$) in the multiscale
system. In this regime the small nonlinear fluctuations of the
multiscale system are relatively large compared to the ramp forcing
for $t$ near zero, so the relative error of the reduced system
responses is large.

\subsection{Predicting the response of the two-scale system via
linear response approximation of the reduced system}
\label{predictingresponse_subsection}

Above in Section \ref{meanresponse_subsection} we discussed the actual
 responses of the statistical mean states of both the two-scale and
reduced models to small Heaviside and ramp forcings. For completeness
of the study, we also attempt to predict the response of the mean
state of the two-scale system via the quasi-Gaussian linear response
approximation
\cite{Abr5,Abr6,Abr7,AbrMaj4,AbrMaj4,AbrMaj5,AbrMaj6,AbrMaj7,MajAbrGro}
of the reduced system. In the quasi-Gaussian response approximation,
the terms $\bar b(\x)$ and $b_A(\x)$ in \eqref{eq:R2} are replaced
with the Gaussian approximations with same mean state and covariance
matrices as in the actual dynamics. This, and the fact that $h(\x)=\x$
in \eqref{eq:h} yields the following formula for the linear response
approximation of the mean state response:
\begin{subequations}
\label{fdt-resp}
\begin{equation}
\begin{array}{c}
\displaystyle
\delta\langle\x\rangle(t)=\int_0^t\R(t-s)\delta\f(s)\ds,\\
\displaystyle
\R(t)=\lim_{r\rightarrow\infty}\frac{1}{r}\int_0^r\x(\tau+t))(\x(\tau)-\bx)^T\dtau\,\Sigma^{-1},
\end{array}
\end{equation}
\begin{equation}
\begin{array}{c}
\displaystyle
\delta\langle\x\rangle_A(t)=\int_0^t\R_A(t-s)\delta\f(s)\ds,\\
\displaystyle
\R_A(t)=\lim_{r\rightarrow\infty}\frac{1}{r}\int_0^r\x_A(\tau+t))(\x_A(\tau)-\bx_A)^T\dtau\,\Sigma_A^{-1},
\end{array}
\end{equation}
\end{subequations}
where $\bx$ and $\Sigma$ are the mean state and covariance matrix of
the corresponding unperturbed systems (two-scale and reduced),
computed as
\begin{subequations}
\begin{equation}
\bx = \lim_{r\rightarrow\infty}\frac{1}{r}\int_0^r\x(\tau)\dtau,
\end{equation}
\begin{equation}
\Sigma = \lim_{r\rightarrow\infty}\frac{1}{r}\int_0^r(\x(\tau)-\bx)
(\x(\tau)-\bx)^T\dtau.
\end{equation}
\end{subequations}
For large multiscale problems the mean response may be difficult to
compute directly, since the large ensemble size needed for an accurate
average is compounded by an already large number of variables and
small timestep discretization. In the case where the mean response of
the slow variables is desired, one might prefer to compute the FDT
response (\ref{fdt-resp}) using a time series from the reduced system
for a ``quick and dirty'' approximation to the ideal response operator
for the multiscale slow system. We show the accuracy of this FDT
response approximation for the Lorenz '96 system, using a
quasi-Gaussian approximation with time series data from the zero- and
first-order reduced systems. The quasi-Gaussian response snapshots for
the response times $T=2$ and $T=5$ are shown in Figures
\ref{fdt_ideal} and \ref{fdt_ideal_ramp} for the Heaviside and ramp
forcing, respectively. Qualitatively, the quasi-Gaussian response does
capture the large features of the actual response, although most noticeable in
these snapshots is the large exaggeration of the quasi-Gaussian
response calculated from the first-order reduced system, which
predicts a much larger off-diagonal response than what is
observed. The possible reason for that is that the distribution
densities of the first-order reduced model are more strongly
non-Gaussian than those of the zero-order reduced model, while the
time autocorrelation functions are more weakly decaying (see Figure
\ref{DDFs}). It was observed previously in \cite{MajAbrGro} that in
these conditions the quasi-Gaussian linear response approximation
tends to overshoot the off-diagonal response by a large margin. In
other words, the better precision of the quasi-Gaussian linear
response of the zero-order model is the result of mutual cancellation
of the two errors: first one is the error in the distribution density
of the zero-order reduced model (significantly more Gaussian than in
in the multiscale dynamics), while the second one is the error in the
quasi-Gaussian linear response due to non-Gaussianity of the
statistical state (less in the case of the zero-order model).

The relative error and cosine similarity are measured against the
multiscale ideal response and can be seen for Heaviside step forcing
in Figure \ref{fdt_ideal_error} and for ramp forcing in Figure
\ref{fdt_ideal_ramp_error}.  The ideal response of the first-order
reduced system is clearly the best of the four responses at capturing
response in the slow variables. It is interesting, but perhaps not too
surprising, that the least accurate estimate is given by the FDT
response of the first-order reduced system. This should be expected
since the quasi-Gaussian approximation is only valid for well-mixing
systems whose distribution densities are close to Gaussian, which in
particular is the case for the uncoupled Lorenz '96 systems in a
chaotic regime $F\geq 8$. However, such a system exhibits suppressed
chaos when coupled to another chaotic systems, and the resulting
distribution density will be far from Gaussian \cite{Abr8}. Since the
first-order reduced system matches more closely the multiscale system,
and the zero-order system will behave as an uncoupled Lorenz '96
system, the first-order system will be less chaotic and will be a poor
candidate for the quasi-Gaussian FDT response. In fact, for
non-chaotic regimes, as in the case of $F_x=7,\ F_y=12$ where
spatially periodic solutions emerge in the two-scale and first-order
reduced systems, the long-time covariance matrix $\Sigma$ will be
singular, so the quasi-Gaussian response as presented will not be applicable. 
For further reading on blended FDT responses which might be
more effective in these cases, see \cite{AbrMaj5}.

\section{Summary}\label{sec:summary}

In this work we studied the response to small external perturbations of
multiscale dynamics and their reduced models for slow variables
only. We elucidated a set of criteria for statistical properties of the
multiscale and reduced systems which facilitated similarity of
responses of both systems to small external perturbations. It was shown
that the similarity of marginal distribution densities of slow
variables and their time autocorrelation functions controlled the
similarity of responses to small external perturbations of both
systems. 

Like in \cite{Abr9}, here we demonstrated that including a first-order
correction term to a standard closure approximation for a nonlinear
chaotic two-time system offered distinct improvements over the
zero-order closure in capturing large-scale features of the slow
dynamics. In particular, this reduced system was able to accurately
capture the distribution density of solutions as well as the mean
state response of the system to simple forcing perturbations. This
correction term was relatively easy to generate, requiring only simple
statistical calculations of the uncoupled fast system for an
appropriate set of fixed parameters, and the resulting reduced system
required much less computational resources than the corresponding
multiscale system.

Focusing on the mean state linear response of the slow variables, we
showed that forcing perturbations in the reduced systems have similar
responses as in the two-time system. Furthermore, we showed that using
the unperturbed dynamics of the reduced systems for linear response
prediction is also possible. However, in the parameter regimes we
present here the first-order reduced systems are not rapidly mixing
and do not follow a Gaussian distribution, but the zero-order reduced
systems do have these properties, so this fluctuation-dissipation
response is effective only using the zero-order system. A linear
response method which takes into account the non-Gaussianity of the
invariant statistical state (such as the blended response algorithm
\cite{AbrMaj4,AbrMaj5,AbrMaj6,AbrMaj7}, based on the tangent map
linear response approximation) is apparently needed to capture the
response for strongly non-Gaussian dynamical regimes in reduced
models.

Here the linear response closure derivation and numerical results have
been presented only for the special case of linear coupling between
slow and fast systems, but this derivation has been extended to
systems with nonlinear and multiplicative coupling \cite{Abr10}. In
future work we hope to extend similar results to these more general
systems and to test the robustness of this method for application to a
large variety of problems.

\medskip

{\bf Acknowledgment.} Rafail Abramov was supported by the National
Science Foundation CAREER grant DMS-0845760, and the Office of Naval
Research grants N00014-09-0083 and 25-74200-F6607. Marc Kjerland was
supported as a Research Assistant through the National Science
Foundation CAREER grant DMS-0845760.

\end{document}